\documentclass[11pt]{article}
\usepackage{latexsym,amssymb,amsmath,amscd,amsthm,amsxtra}
\usepackage{mathrsfs}
\usepackage{epsfig}
\usepackage{pgf,tikz}
\usepackage{graphicx}
\usepackage{graphics}
\usepackage{setspace}
\newcommand{\D}{\displaystyle}
\newcommand{\counte}{section}
\newtheorem{define}{\bf Definition}[\counte]

\newtheorem{prop}{\bf Proposition}[\counte]

\newtheorem{theorem}{\bf Theorem}[\counte]
\newtheorem{coro}{\bf Corollary}[\counte]

\newtheorem{remark}{\bf Remark}[\counte]

\author{Zhao Xu-an, zhaoxa@bnu.edu.cn\\Gao Hongzhu, hzgao@bnu.edu.cn\\
Department of Mathematics, Beijing Normal University\\Key Laboratory
of Mathematics and Complex Systems\\ Ministry of Education,
China, Beijing 100875}

\title{Differential geometry of general affine plane curves\thanks{The authors are supported by NSFC 11571038}}
\date{}

\begin{document}

\maketitle
\begin{abstract}
In this paper we study the general affine geometry of curves in affine space $A^2$. For a regular plane curves we define two kinds of moving frames. The first is of minimal order in all moving frames.
The second is the Frenet moving frame. We get the moving equations of these moving frames. And we prove that curvature and signature are the complete invariants of regular curves. As application we give
a complete classification of constant curvature curves in $A^2$. %We also give compute the similar derivatives of the equation of the curve. on and compute the
\end{abstract}

\noindent{\bf Keywords: }General affine differential geometry, Plane curve, Moving frame, Invariant arc element, Curvature.

\noindent{\bf MSC(2010): }Primary 53B52

\section{Introduction}

The affine geometry was founded by Blaschke, Pick, Radon, Berwald and Thomsen among others in the period
from 1916 to 1923. And a systematic theory of curves and surfaces in three dimension affine space was developed. For accounts and expository books
appeared on the subject, see Blaschke\cite{Blaschke_85}, Guggenheimer\cite{Guggenheimer_63} and Spivak's\cite{Spivak}. In these works the affine geometry means the equiaffine geometry. That is the Kleinnian geometry of
the group of affine transformations which preserve volume. But for general affine geometry, that is the Kleinian geometry of the group of general affine transformations, there is little work. All the work we can find are Weise\cite{Weise_38}\cite{Weise_39}, Kllingenberg\cite{Klingenberg_81_1}\cite{Klingenberg_51_2}, Svec\cite{Svec_59},
Wilkinson\cite{Wilkinson_88} and Weiner\cite{Weiner_94}. Weise and Klingenberg began the study of general affine differential geometry. %, but their works did not continue by other mathematicians
Wilkinson discussed submanifolds of low codimension. Svec and Weiner studied surfaces in affine space $A^4$ independently. All the discussions in these papers were in an abstract form. %The study of general affine geometry still has a long way to go.
We need a theory of general affine geometry which is parallel to the classical Euclidean geometry of curves and surfaces. And we haven't found any reference on the explicit computation even for curves in $A^2$. This is the motivation we write this paper.

The content of this paper is as follows. In section 2 we give an introduction to jet spaces and Fels and Olver's moving frame method. In section 3 we construct a moving frame of minimal order. In section 4 we compute the arc element and the curvature of a regular curve and get the moving equations of the moving frame. In section 5 we construct the Frenet moving frame. In section 6 we give an affine classification curves of constant curvature. In section 7 we discuss modular invariants.

\section{Jet spaces and the moving frame method}

The local differential geometry of a curve $C$ at a point $P$ is determined by the shape of $C$ at an arbitrary small neighborhood. So all the local differential geometric properties and invariants of $C$ are determined by the local
data of $C$. It is useful to isolate the informations of curve $C$ at $P$. This idea hints the concept of jet of curve. The general definition of jet of submanifolds was given by Ehresmann at 1950s for an ambient manifold $M$. %See Ehresmann\cite{Olver_95} for reference.

For a smooth manifold $M$, let $SM^{d}_P(M)$ be the set of all smooth $d$ dimensional submanifolds of $M$ that contain the point $P$. For integer $r\geq 0$, define an equivalence relation $\sim$ on $SM^{d}_P(M)$ such that $N_1\sim N_2$ iff $N_1$ and $N_2$ have contact at least of order $r$. The Jet space of $d$-dimensional submanifolds of $M$ at the point $P$ of order $r$ is the quotient set $J_P^{d,r}(M)=SM^{d}_P(M)/\sim$. And $J^{d,r}(M)=\bigcup\limits_{P\in M} J_P^{d,r}(M)$ is the jet space of $d$-dimensional submanifolds of $M$ of order $r$. An element of jet space $J^{d,r}(M)$ is called a $d$-jet of order $r$. For an $n$ dimensional manifold $M$, $J^{d,r}(M)$ is also a manifold, and the dimension of jet space $J^{d,r}(M)$ is $\displaystyle{d+(n-d){d+r\choose r}}$. In this paper we use the Jet space $J^{1,r}(A^2)$ for $r>0$ in the study of curves in $A^2$.

\subsection{Action of $Aff(2)$ on $J^{1,r}(A^2)$}
Let $Aff(2)$ be the group of general affine transformations of affine plane $A^2$. %The dimension of $Aff(2)$ is $6$.
The general affine geometry is the Kleinian geometry given by the group action $Aff(2)\times A^2\to A^2$. Under the affine coordinates $(x,y)$ on $A^2$, a general affine transformation has the form
$$\left(
\begin{array}{c}
 x' \\
  y' \\
   \end{array}
   \right)=\left(
             \begin{array}{cc}
               a_{11} & a_{12} \\
               a_{21} & a_{22} \\
             \end{array}
           \right)
   \left(
\begin{array}{c}
 x \\
  y \\
   \end{array}
   \right)+\left(
\begin{array}{c}
 x_0 \\
  y_0\\
   \end{array}
   \right).$$
   Or
$$\left(
\begin{array}{c}
 x' \\
  y' \\ 1 \\
   \end{array}
   \right)=\left(
             \begin{array}{ccc}
               a_{11} & a_{12} & x_0\\
               a_{21} & a_{22} & y_0\\
               0&0&1\\
             \end{array}
           \right)
   \left(
\begin{array}{c}
 x \\
  y \\ 1\\
   \end{array}
   \right).$$

The coordinates on $Aff(2)$ are $a_{11},a_{12},a_{21},a_{22},x_0,y_0$.

In this paper, we assume a curve has the form $y=f(x)$ locally. The local coordinates of the manifold $J^{1,r}(M)$ are $x,y,y_x,y_{xx},\cdots,y_{xx\cdots x}$, where $x$ appears $r$ times in the subindex of the last item. The action of $Aff(2)$ on $A^2$ induces an action on $SM^{1}(A^2)$ and hence on $J^{1,r}(A^2),r\geq 0$.

\subsection{Differential invariants of jet spaces}

By using jet space as a tool, we can transform the study of local affine congruence invariants of $C$ at $P$ to the invariants of jet space under the action of $Aff(2)$. We have the following definition.

\begin{define}
A smooth function $f:J^{1,r}(A^2)\to \mathbb{R}$ which is invariant under the action of $Aff(2)$ is called a differential invariant of curves in $A^2$ of order less than $r+1$.
\end{define}

Jet space can be regarded as the finite dimension cut-off of infinite dimension space of all submanifolds. The use of jet space separates the study of differential geometry into algebra part and analysis part. So it makes the structure of differential geometry theory clearer.

\subsection{Moving frames for jet spaces of curves}

The basic language of modern differential geometry is Cartan's moving frame. For a curve $C$, the moving frame method gives an affine frame at each point $P\in C$.
That is the point $P\in A^2$ and two linearly independent vectors $e_1$ and $e_2$ in the associated vector space of $A^2$. %The data of an affine frame is $\in (P,e_1,e_2)$.
The moving equations of moving frame are given by
$$dP= w_1 e_1+w_2 e_2 $$
$$d\left(
\begin{array}{c}
 e_1 \\
  e_2 \\
   \end{array}
   \right)=\left(
             \begin{array}{cc}
               w_{11} & w_{12} \\
               w_{21} & w_{22} \\
             \end{array}
           \right)
   \left(
\begin{array}{c}
 e_1 \\
  e_2 \\
   \end{array}
   \right).$$
By Cartan's result, if we construct a moving frame on curve $C$ and compute the one forms $w_1,w_2,w_{11},w_{12}$, $w_{21},w_{22}$, then we determine the congruence class of the curve $C$.

\subsection{Fels and Olver's moving frame method}
This section is an introduction to Fels and Olver's moving frame method. See Fels and Olver\cite{Fels_Olver_98}\cite{Fels_Olver_99}\cite{Fels_Olver_04} for details.

\begin{define} Let $\phi:G\times M\to M$ be a smooth action of Lie group $G$ on smooth manifold $M$. A moving frame on $M$ is a smooth, G-equivariant map $\rho: M\to G$.
\end{define}

There are two types of moving frames.

$\left\{
  \begin{array}{ll}
    \rho(gz) =g\rho(z), & \hbox{left moving frame;} \\
    \rho(gz) =\rho(z)g^{-1}, & \hbox{right moving frame.}
  \end{array}
\right.$

\begin{theorem}
A moving frame exists on $M$ if and
only if $G$ acts freely and regularly on $M$.
\end{theorem}

The explicit construction of a moving frame is based on Cartan's normalization procedure.
Let $G$ act freely and regularly on $M$ and $K$ be a cross-section to the group orbits, that is a submanifold $K$ which transversally intersects each orbit once. Let $g$ be the unique group element which maps $P$ into the
cross-section $K$, then $\rho:M \to G,P\mapsto g$ is a right moving frame. And $\rho: M \to G,P\mapsto g^{-1}$ is a left moving frame.
The unique intersection point of the orbit of $P$ and $K$ can be viewed as the canonical form or normal form of $P$, as
prescribed by the cross-section $K$.

If a moving frame is in hand, the determination of the invariants is routine.
The specification of a moving frame by choosing a cross-section induces a canonical
procedure to map functions to invariants.

\begin{define} The invariantization of a function $F: M\to\mathbb{R}$ is the unique invariant
function $\iota(F)$ that coincides with $F$ on the cross-section, that is $\iota(F)|_K= F|_K$.\end{define}

Invariantization defines a projection from the space of (smooth) functions to the space of invariants that, moreover,
preserves all algebraic operations.

%Since the action $Aff(2)\times J^{1,r}(A^2)\to J^{1,r}(A^2)$ is not free, so the moving frame can only be constructed on certain open set of $J^{1,r}(A^2)$. We will discuss this later.

\section{The construction of the moving frame}
By choosing suitable affine coordinates we can write a curve $C$ as $y=y(x)$ locally. %Since any jet can be transformed by elements in $Aff(2)$ into this form, it gives no essential restriction.

Given a jet in $J^{1,r}(A^2)$, we use a transportation to move this jet to a jet at $O(0,0)$. The element of $Aff(2)$ fixing the point $(0,0)$ has the form of $x'=a_{11} x+a_{12}y, y'=a_{21}x+a_{22} y$.
%To construct left moving frame on jet space, we use the transformation $A^{-1}$. The explicit formula of $A^{-1}$ is given by
%$$x'=\frac{a_{22} x-a_{21} y-a_{22} a_1 +a_{21} a_2}{\Delta}, y'=\frac{-a_{12} x+ a_{11} y+a_{12} a_1-a_{11} a_2}{\Delta}$$ Where

By direct computation, we get the action of $Aff(2)$ on $J^{1,5}(A^2)$. The explicit formulae are given by

\begin{equation}\label{y'_1}
\displaystyle{y'_{x'}=\frac{a_{21}+a_{22}y_x} {a_{11}+a_{12}y_x}}
\end{equation}

\begin{equation}\label{y'_2}
\displaystyle{y'_{x'x'}=\frac{\Delta y_{xx}}{\Gamma^3}}
\end{equation}

\begin{equation}\label{y'_3}
\displaystyle{y'_{x'x'x'}=\frac{\Delta (a_{11}y_{xxx}+a_{12}y_x y_{xxx}-3a_{12}y^2_{xx})}{\Gamma^5}}
\end{equation}

\begin{equation}\label{y'_4}
\displaystyle{y'_{x'x'x'x'}=\frac{\Delta M}{\Gamma^7}}
\end{equation}

\begin{equation}\label{y'_5}
\displaystyle{y'_{x'x'x'x'x'}=\frac{\Delta N}{\Gamma^9}}
\end{equation}

\noindent Where $\Delta=a_{11}a_{22}-a_{21}a_{12}$, $\Gamma=a_{11}+a_{12} y_x$,

%$M=a^2_{11}y_{xxxx}+2a_{11}a_{12}y_x y_{xxxx}+a^2_{12}y^2_x y_{xxxx}-10 a_{11}a_{12}y_{xx} y_{xxx}-10 a^2_{12}y_x y_{xx} y_{xxx} +15a^2_{12}y^3_{xx}$,

$M=\Gamma^2 y_{xxxx}
-10a_{12} \Gamma y_{xx} y_{xxx}
+15a_{12}^2 \Gamma y_{xx}^3,$

%$N=a^3_{11}y_{xxxxx}+3a^2_{11} a_{12} y_x y_{xxxxx}+3a_{11}a^2_{12}y^2_x y_{xxxxx}+a^3_{12} y^3_x y_{xxxxx}-15a^2_{11} a_{12} y_{xx} y_{xxxx}$

$N=\Gamma^3 y_{xxxxx}
-15a_{12} \Gamma^2 y_{xx}y_{xxxx}
-10a_{12} \Gamma^2 y_{xxx}^2
+105a_{12}^2 \Gamma y_{xx}^2 y_{xxx}
-105a_{12}^3 y_{xx}^4.$

%$-30a_{11} a^2_{12}  y_x y_{xx} y_{xxxx}-15a^3_{12} y^2_x y_{xx} y_{xxxx}-10a^2_{11} a_{12} y^2_{xxx}-20a_{11}a^2_{12} y_x y^2_{xxx}-10a^3_{12} y^2_x y^2_{xxx}$

%$+105a_{11}a^2_{12} y^2_{xx} y_{xxx} +105a_{12}^3 y_x y^2_{xx} y_{xxx}-105a^3_{12} y^4_{xx}$

%These formulae are derived from the chain rule of composition of maps and give the action of $Aff(2)$ on the jet coordinates of manifold $J^{1,5}(A^2)$.

\subsection{Explicit construction of moving frame}

In this section we use Fels and Olver's moving frame method to construct a right moving frame on jet space $J^{1,4}(A^2)$ and simultaneously we show how to transform a jet of curve at $P(x,y)$ into a standard jet at $(0,0)$.

1. Since $\displaystyle{y'_{x'}=\frac{a_{21}+a_{22}y_x} {a_{11}+a_{12}y_x}}$, if we fix $y'_{x'}=0$, then we get
\begin{equation}\label{a22a12}
a_{21}+a_{22}y_x=0
\end{equation} and $a_{11}+a_{12}y_{x}\not=0$. So we have
\begin{equation}\label{y_1}
\displaystyle{y_x=-\frac{a_{21}}{a_{22}}}
\end{equation}
Here we require $a_{22}\not=0$. Otherwise $a_{21}=-a_{22}y_x=0$, as a consequence $\Delta=0$. This is a contradiction.

Under this condition we have $\displaystyle{\Gamma=a_{11}+a_{12}y_x=a_{11}+a_{12}\frac{-a_{21}}{a_{22}}=\frac{\Delta}{a_{22}}}$.

2. Since $\displaystyle{y'_{x'x'}=\frac{\Delta y_{xx}}{\Gamma^3}}$, if we assume $y_{xx}\not=0$, then it is reasonable to set $y'_{x'x'}=1$. Substituting $\displaystyle{\Gamma=\frac{\Delta}{a_{22}}}$ into $\displaystyle{y'_{x'x'}=\frac{\Delta y_{xx}}{\Gamma^3}}$, we get
\begin{equation}\label{a22delta}
y_{xx}a_{22}^3=\Delta^2.
\end{equation}
Hence
\begin{equation}\label{y_2}
y_{xx}=\frac{\Delta^2} {a_{22}^3}.
\end{equation}

This require that $y_{xx}$ has the same sign with $a_{22}$.

%When we fix $y'_{x'}$ and $y'_{x'x'}$, then $y_x$ and $y_{xx}$ satisfy the relation $$\displaystyle{y_x=\frac{a_{12}}{a_{11}}}.\label{y_1}$$ $$y_{xx}=\frac{\Delta^2}{a_{11}^3}. \label{y_2}$$

3. We substitute the expressions of $y_x$ and $y_{xx}$ in Equation \ref{y_1} and \ref{y_2} into Equation \ref{y'_3} and get
$$\displaystyle{y'_{x'x'x'}=a_{22}^4\Delta^{-3}(y_{xxx}-\frac{3a_{12}\Delta^3}{a^5_{22}})}.$$

We fix $y'_{x'x'x'}=0$. That is
\begin{equation}\label{y_3}
y_{xxx}=\frac{3a_{12}\Delta^3}{a^5_{22}}
\end{equation} Hence

\begin{equation}\label{a12delta}
a_{12}=\frac{1}{3}a^5_{22}\Delta^{-3}y_{xxx}
\end{equation}

4. As in 3, we substitute the expressions of $y_x,y_{xx}$ and $y_{xxx}$ in Equation \ref{y_1}, \ref{y_2} and \ref{y_3} into Equation \ref{y'_4} and get
$$\displaystyle{y'_{x'x'x'x'}=a_{22}^5 \Delta^{-4}(y_{xxxx}-\frac{15a^2_{12}\Delta^4}{a^7_{22}})}=\frac{1}{a_{22} y^2_{xx}}(y_{xxxx}-\frac{5y^2_{xxx}}{3y_{xx}})$$

Since we have assumed $y_{xx}\not=0$, and $y_{xx}$ has the same sign with $a_{22}$, if we further assume $\displaystyle{y_{xxxx}-\frac{5y^2_{xxx}}{3y_{xx}}\not=0}$, $y'_{x'x'x'x'}$ must have the same sign with $\displaystyle{y_{xx}(y_{xxxx}-\frac{5y^2_{xxx}}{3y_{xx}})}$ $\displaystyle{=\frac{1}{3}(3y_{xx}y_{xxxx}-5y^2_{xxx})}$.

We fix $y'_{x'x'x'x'}=\left\{
             \begin{array}{ll}
               1, & \hbox{if} \ \ 3y_{xx}y_{xxxx}-5y^2_{xxx}>0; \\
               -1, & \hbox{if}  \ \ 3y_{xx}y_{xxxx}-5y^2_{xxx}<0.
             \end{array}
           \right. $

In either cases we have
\begin{equation}\label{a11}
\displaystyle{a_{22}=\frac{1}{ 3y^3_{xx}}|3y_{xx}y_{xxxx}-5y^2_{xxx}|}
\end{equation}

%For this to be true, we must have $\displaystyle{y_{xxxx}-\frac{5y^2_{xxx}}{3y_{xx}}\not=0}$.
In the previous four steps we transform the jet of curve $C$ at $P(x,y)$ into a jet with $x'=y'=y'_{x'}=y'_{x'x'x'}=0$, $y'_{x'x'}=1$ and $y'_{x'x'x'x}=\pm 1$ if both $y_{xx}$ and $3y_{xx}y_{xxxx}-5y^2_{xxx}$ are not $0$. In the following we compute the affine transformation to get the moving frame.

From Equation \ref{a22a12}, we have

\begin{equation}\label{a21}
\displaystyle{a_{21}=-\frac{y_x}{ 3y^3_{xx}}|3y_{xx} y_{xxxx}-5y^2_{xxx}|}
\end{equation}

If we assume $\Delta>0$, then from Equation \ref{a22delta}  and 12 we have

\begin{equation}\label{delta}
\displaystyle{\Delta=\sqrt{y_{xx}a_{22}^3}=\sqrt{\frac{1}{ 27y^8_{xx}}|3y_{xx}y_{xxxx}-5y^2_{xxx}|^3}}
\end{equation}

%$\displaystyle{a_{21}=\frac{y_{xxx}}{3y^{10}_{xx}}(y_{xxxx}-\frac{5y^2_{xxx}}{3y_{xx}})^5 /\sqrt{\frac{1}{ y^{15}_{xx}}(y_{xxxx}-\frac{5y^2_{xxx}}{3y_{xx}})^9}}=\frac{y_{xxx}}{3y^{2}_{xx}}(y_{xxxx}-\frac{5y^2_{xxx}}{3y_{xx}}) /\sqrt{y_{xx}(y_{xxxx}-\frac{5y^2_{xxx}}{3y_{xx}})}$

Substituting $\Delta$ into Equation \ref{a12delta}, we get
\begin{equation}\label{a12}
\displaystyle{a_{12}=\frac{1}{3\sqrt{3}}\frac{y_{xxx}|3y_{xx} y_{xxxx}-5y^2_{xxx}|^{\frac{1}{2}}}{y^3_{xx}}}
\end{equation}

By virtue of $a_{11}a_{22}-a_{21}a_{12}=\Delta$, we have
%If $a_{11}$ and $y_{xx}$ have different symbol, then we have $3 y_{xx}y_{xxxx}-5y_{xxx}^2<0$, then we should require $y'_{x'x'x'x'}=-1$.
\begin{equation}\label{a22}
\displaystyle{a_{11}=-\frac{1}{3\sqrt{3}}\frac{(y_x y_{xxx}-3y^2_{xx})|3y_{xx} y_{xxxx}-5y^2_{xxx}|^{\frac{1}{2}}}{y^3_{xx}}}
\end{equation}
%¶ÔÓÚ$y_{xx}=0$µÄÇéÐΣ¬×ÔÈ»$y'_{x'x'}=0$.´Ëʱ

%$\displaystyle{y'_{x'x'x'}=\frac{\Delta (a_{11}y_{xxx}+a_{21}y_x y_{xxx}-3a_{21}y^2_{xx})}{\Gamma^5}}=a_{22}^4\Delta^{-3}y_{xxx}$.

%$\frac{\sqrt{3}(27y_2^4 y_4^2 +18y_1 y_2^2 y_3 y_4^2 -90y_2^3 y_3^2 y_4 -60y_1 y_2 y_3^3 y_4 +75y_2^2 y_3^4 +50y_1 y_3^5)}{{\sqrt {-{\frac { ( 3y_2 y_4-5y_3^2)^{3}}{y_2^8}}}} y_2^7}$

%In any cases we have

%$\displaystyle{a_{11}=\frac{1}{ y^3_{xx}}|y_{xx} y_{xxxx}-\frac{5}{3}y^2_{xxx}|}$

%$\displaystyle{a_{21}=-\frac{1}{3}\frac{y_{xxx}|y_{xx} y_{xxxx}-\frac{5}{3}y^2_{xxx}|^{\frac{1}{2}}}{y^3_{xx}}}$

%$\displaystyle{a_{12}=\frac{y_x}{ y^3_{xx}}|y_{xx} y_{xxxx}-\frac{5}{3}y^2_{xxx}|}$

%$\displaystyle{a_{22}=-\frac{1}{3}\frac{(y_x y_{xxx}-3y^2_{xx})|y_{xx} y_{xxxx}-\frac{5}{3}y^2_{xxx}|^{\frac{1}{2}}}{y^3_{xx}}}$
Summarizing the results in the previous discussions, and setting $S_1=y_{xx}$, $S_2=3y_{xx} y_{xxxx}-5y^2_{xxx}$, then we have

\begin{prop}
Let $j$ be a jet of curve $y=y(x)$ at $(x_0,y_0)$ with $S_1\not=0$ and $S_2\not=0$, then there exists a unique element $A\in Aff(2)$ with $\Delta>0$ such that $Aj$ is a jet of curve $y'=y'(x')$ satisfying $x'=y'=y'_{x'}=y'_{x'x'x'}=0, y'_{x'x'}=1$ and $y'_{x'x'x'x'}=\pm 1$. Furthermore $A$ has the form
$$\left(
\begin{array}{c}
 x' \\
  y' \\
   \end{array}
   \right)=\left(
             \begin{array}{cc}
               a_{11} & a_{12} \\
               a_{21} & a_{22} \\
             \end{array}
           \right)
   \left(
\begin{array}{c}
 x-x_0 \\
  y-y_0 \\
   \end{array}
   \right)$$
Where $a_{22},a_{21},a_{12},a_{11}$ are given by Equation \ref{a11},\ref{a21},\ref{a12} and \ref{a22}.
\end{prop}

By this proposition, we give the following definition.

\begin{define}
1. Let $C:y=y(x)$ be a smooth curve. A point $P(x,y)\in C$ satisfying $S_1\not=0$ and $S_2\not=0$ is called a regular point of $C$. A point which is not regular is called a
singular point.

2. A jet $j\in J^{1,4}(A^2)$ is regular if there exist affine coordinates $x,y$ such that $j$ has the form $y=y(x)$ locally with both $S_1\not=0$ and $S_2\not=0$. A jet $j\in J^{1,4}(A^2)$ is singular if it is not regular.
\end{define}

The definition of regularity of a jet does not depend on the affine coordinates. Direct computation shows under affine transformations, we have $S'_1=\Delta\Gamma^{-3} S_1$, $S'_2=\Delta^2\Gamma^{-8}S_2$.

%$\displaystyle{y'_{x'x'}=\Delta\Gamma^{-3}}y_{xx}$;

%$3y'_{x'x'}y'_{x'x'x'x'}-5y'^2_{x'x'x'}=\Delta^2\Gamma^{-8}(3y_{xx}y_{xxxx}-5y^2_{xxx}).$

%\begin{coro}A jet \end{coro}

%We Choosing the moving frame $P,e_1,e_2$ with

%$e_1=\left(
     %      \begin{array}{c}
      %       \displaystyle{\frac{1}{ 3y^3_{xx}}|3y_{xx}y_{xxxx}-5y^2_{xxx}|}\\
       %      \displaystyle{\frac{y_x}{ 3y^3_{xx}}|3y_{xx} y_{xxxx}-5y^2_{xxx}|} \\
      %     \end{array}
      %   \right)$

%$e_2=\left(
       %    \begin{array}{c}
        %     \displaystyle{-\frac{1}{3\sqrt{3}}\frac{y_{xxx}|3y_{xx} y_{xxxx}-5y^2_{xxx}|^{\frac{1}{2}}}{y^3_{xx}}} \\
        %     \displaystyle{-\frac{1}{3\sqrt{3}}\frac{(y_x y_{xxx}-3y^2_{xx})|3y_{xx} y_{xxxx}-5y^2_{xxx}|^{\frac{1}{2}}}{y^3_{xx}}} \\
       %    \end{array}
      %   \right)$

\begin{prop} If a smooth curve $C:y=y(x)$ is regular, then we have a smooth right moving frame $(P,\alpha_1,\alpha_2)$ on $C$,
where

$\left(
         \begin{array}{c}
           \alpha_1 \\
           \alpha_2 \\
         \end{array}
       \right)=
\left(
             \begin{array}{cc}
               a_{11} & a_{12} \\
               a_{21} & a_{22} \\
             \end{array}
           \right)=\left(
             \begin{array}{cc}
               \displaystyle{\frac{-\sqrt{3}(y_x y_{xxx}-3y^2_{xx})|3y_{xx} y_{xxxx}-5y^2_{xxx}|^{\frac{1}{2}}}{9y^3_{xx}}} & \displaystyle{\frac{\sqrt{3}y_{xxx}|3y_{xx} y_{xxxx}-5y^2_{xxx}|^{\frac{1}{2}}}{9y^3_{xx}}}  \\
               -\displaystyle{\frac{y_x}{ 3y^3_{xx}}|3y_{xx} y_{xxxx}-5y^2_{xxx}|} &  \displaystyle{\frac{1}{ 3y^3_{xx}}|3y_{xx}y_{xxxx}-5y^2_{xxx}|} \\
             \end{array}
           \right)$
\end{prop}

%The left moving frames is given by the inverse the right moving frames.
%\section{Plane affine differential geometry}

\begin{prop}
If all the points of curve $C$ is singular, then there exists $C_0,C_1,C_2,C_3$ such that $y=(C_3 x+C_2)^{\frac{1}{2}}+C_1 x +C_0$ or there exists $C_0,C_1,C_2$ such that $y=C_2 x^2+ C_1 x+C_0$.
Equivalently $C$ is affine congruence to parabola $y=x^2$ or line $y=x$.
\end{prop}

\noindent Proof: The singular condition is $ \displaystyle{3y_{xxxx}y_{xx}-5y^2_{xxx}}=0$. It can be written as $\displaystyle{\frac{3y_{xxxx}}{y_{xxx}}=\frac{5y_{xxx}}{y_{xx}}}$, if $y_{xxx}\not=0$, we have $3\ln|y_{xxx}|=5\ln|y_{xx}|$. Solving this differential equation will finish the proof.

\section{Invariant arc element, curvature and the moving equations of moving frame}

In this section, we derive the invariant arc element and curvature for an affine plane curve. %We use the Cartan's normalization method.

\begin{prop}
The affine invariant arc element is $ds=\displaystyle{\frac{|3y_{xx}y_{xxxx}-5y^2_{xxx}|^{\frac{1}{2}}}{\sqrt{3}y_{xx}}dx}$. As a consequence $\displaystyle{\frac{d}{ds}=\frac{\sqrt{3}y_{xx}}{|3y_{xx}y_{xxxx}-5y^2_{xxx}|^{\frac{1}{2}}}\frac{d}{dx}}$ is an affine invariant differential operator.
\end{prop}

\noindent {\bf Proof: }We use the right moving frames to compute the affine invariant arc element. Differentiate $x'= a_{11}x + a_{12}y$, we get $dx'=a_{11}dx+a_{12}dy$. On the curve $C$, we have $dy=y_{x}dx$. By the Normalization procedure, we substitute $dy$ and $a_{11},a_{21}$ into the expression of $dx'$ and get $dx'=\displaystyle{\frac{|3y_{xx}y_{xxxx}-5y^2_{xxx}|^{\frac{1}{2}}}{\sqrt{3}y_{xx}}}dx$.
Since $dx'$ is the invariantization of $dx$, it can be chosen as the invariant arc element. Dually $\displaystyle{\frac{d}{dx'}}$ is an invariant differential operator.

\begin{theorem}
The curvature of curve $C: y=y(x)$ at the point $P(x,y)$ is
$$\displaystyle{k(x)=\frac{\sqrt{3}(9y_{xx}^2 y_{xxxxx}-45y_{xx}y_{xxx}y_{xxxx}+40y_{xxx}^3)}{3|3y_{xx}y_{xxxx}-5y^2_{xxx}|^{\frac{3}{2}}}}.$$
\end{theorem}

\noindent {\bf Proof: }Substituting Equations \ref{a11},\ref{a21},\ref{a12},\ref{a22} into Equation \ref{y'_5}, we get the invariantization of $y_{xxxxx}$ which is the curvature.

\begin{theorem}
The moving equation of the right moving frame is $$\displaystyle{\frac{d}{ds} \left(
\begin{array}{c}
 \alpha_1 \\
  \alpha_2 \\
   \end{array}
   \right)=\left(
                       \begin{array}{cc}
                         \frac{1}{2}k(s) & \frac{1}{3} \sigma \\
                         -1 & k(s) \\
                       \end{array}
                     \right)\left(
\begin{array}{c}
 \alpha_1 \\
  \alpha_2 \\
   \end{array}
   \right)}
$$ Where $\sigma=sgn(S_2),\alpha_1=(\displaystyle{\frac{-\sqrt{3}(y_x y_{xxx}-3y^2_{xx})|3y_{xx} y_{xxxx}-5y^2_{xxx}|^{\frac{1}{2}}}{9y^3_{xx}}}, \displaystyle{\frac{\sqrt{3}y_{xxx}|3y_{xx} y_{xxxx}-5y^2_{xxx}|^{\frac{1}{2}}}{9y^3_{xx}}}) $ and \\ $\alpha_2=(-\displaystyle{\frac{y_x}{ 3y^3_{xx}}|3y_{xx} y_{xxxx}-5y^2_{xxx}|}, \displaystyle{\frac{1}{ 3y^3_{xx}}|3y_{xx}y_{xxxx}-5y^2_{xxx}|})$
%Where $\displaystyle{k(x)=\pm\frac{\sqrt{3}(9y_{xx}^2 y_{xxxxx}-45y_{xx}y_{xxx}y_{xxxx}+40y_{xxx}^3)}{3|3y_{xx}y_{xxxx}-5y^2_{xxx}|^{\frac{3}{2}}}}$.
\end{theorem}

\noindent {\bf Proof: }This is derived by direct computation from the results in Proposition 3.2, Proposition 4.1 and Theorem 4.1.

The left moving frame is given by $\displaystyle{(e_1,e_2)=\left(
                                                                          \begin{array}{cc}
                                                                            a_{11} & a_{12} \\
                                                                            a_{21} & a_{22} \\
                                                                          \end{array}
                                                                        \right)^{-1}}$.

\begin{theorem}
The moving equation of the left moving frame is $$\displaystyle{\frac{d}{ds} (e_1,e_2)=(e_1,e_2) \left(
                       \begin{array}{cc}
                         -\frac{1}{2}k(s) & -\frac{1}{3} \sigma \\
                         1 &- k(s) \\
                       \end{array}
                     \right)}
$$ Where $e_1=\left(
\begin{array}{c}
 \displaystyle{\frac{\sqrt{3}y_{xx}}{|3y_{xx} y_{xxxx}-5y^2_{xxx}|^{\frac{1}{2}}}} \\
  \displaystyle{\frac{\sqrt{3}y_x y_{xx}}{|3y_{xx} y_{xxxx}-5y^2_{xxx}|^{\frac{1}{2}}}} \\
   \end{array}   \right)$, $e_2=\left(
\begin{array}{c}
 \displaystyle{-\frac{y_{xx}y_{xxx}}{|3y_{xx} y_{xxxx}-5y^2_{xxx}|}} \\
   \displaystyle{-\frac{y_{xx}(y_x y_{xxx}-3y^2_{xx})}{|3y_{xx} y_{xxxx}-5y^2_{xxx}|}} \\
   \end{array}   \right)$
\end{theorem}

\noindent {\bf Proof: }Let $A=\left(
                                                                          \begin{array}{cc}
                                                                            a_{11} & a_{12} \\
                                                                            a_{21} & a_{22} \\
                                                                          \end{array}
                                                                        \right)$ and $K=\left(
                       \begin{array}{cc}
                         \frac{1}{2}k(s) & \frac{1}{3}\sigma \\
                         -1 &k(s) \\
                       \end{array}
                     \right)$. We have $dA^{-1}=-A^{-1}dA A^{-1}=-A^{-1}KA A^{-1}=-A^{-1}K$.

By Cartan's moving frame method, Theorem 4.2 or Theorem 4.3 gives the general affine classification of regular plane curves.

\begin{coro}
For a regular smooth plane curve, $k=k(s)$ and $\sigma$ are the complete affine invariants. All the differential invariants are functions of $k(s)$ and its derivatives on $s$.
\end{coro}

\begin{remark}
The discuss in this section is for the geometry of orientation preserving affine transformations. %For the general affine transformation group $Aff(2)$,
If a transformation reverses the orientation, then it changes the sign of $k(s)$ and keeps $\sigma$ invariant.
\end{remark}

\section{The Frenet moving frame}
The left moving frame $e_1,e_2$ we constructed in the last section uses the derivatives of $y(x)$ up to order $4$ and the moving equations of the moving frame contains derivatives up to order $5$.

Using arc element we have a left moving frame which contains the derivatives of $y(x)$ up to order $5$, so the moving equations of moving frame contains derivatives of order $6$. We call it Frenet frame. Computation gives the following result.

\begin{theorem}
Let $r=\left(
                                 \begin{array}{c}
                                   x \\
                                   y(x) \\
                                 \end{array}
                               \right)$, $\displaystyle{t=\frac{dr}{ds}}$ and $\displaystyle{n=\frac{dt}{ds}}$, then

$$\frac{d}{ds}(t,n,r)=(t,n,r)\left(
                       \begin{array}{ccc}
                         0 & -\frac{1}{2} k^2-\frac{1}{2}\frac{dk}{ds}-\frac{\sigma}{3} &  1\\
                         1 & -\frac{3}{2}k& 0\\   0&0&0\\
                       \end{array}
                     \right)
$$
\end{theorem}

\noindent {\bf Proof: }Let $\displaystyle{r=\left(
                                 \begin{array}{c}
                                   x \\
                                   y(x) \\
                                 \end{array}
                               \right)}$,
                               computation shows $\displaystyle{t=\frac{dr}{ds}=e_1}$. We set $\displaystyle{\frac{dt}{ds}=n}$. So $\displaystyle{n=\frac{de_1}{ds}=-\frac{1}{2}k(s)t}$ $+e_2$ and
$\displaystyle{\frac{dn}{ds}=(-\frac{1}{2} k^2-\frac{1}{2}\frac{dk}{ds}-\frac{\sigma}{3})t-\frac{3}{2}k n}$.

The moving frame $e_1,e_2$ is valid for curve with derivative up to order $5$, but the Frenet frame $t,n$ is only valid for curve with derivative up to order $6$.

\begin{define}
The vectors $t$ and $n$ are respectively called the affine tangential vector and normal vector of the curve $C$.
\end{define}

%\begin{remark}\end{remark}

%\begin{remark}\end{remark}

\section{Curves with constant curvature}
In this section we give the classification of curves with constant curvature in $A^2$. We assume all curves are regular. The one parameter subgroup of $Aff(2)$ has the form $\exp(tX)$, where $X\in aff(2)$ and $aff(2)$ is the Lie algebra of $Aff(2)$. We write $A\in Aff(2)$ as a $3\times 3$ matrix $\left(
                                                                                           \begin{array}{ccc}
                                                                                             a_{11} & a_{12} & x_0 \\
                                                                                             a_{21} & a_{22} & y_0 \\
                                                                                             0 & 0 & 1 \\
                                                                                           \end{array}
                                                                                         \right)$.

The following proposition gives the one parameter subgroups of $Aff(2)$ up to conjugation.

\begin{prop} Let $G=Exp(tX)$ be a one parameter subgroup of $Aff(2)$, %then $X$ has the following form.

1. If $G$ fixes a point $P\in A^2$, then up to conjugation $X$ has the form of

$\left(
               \begin{array}{ccc}
                 \lambda_1 & 0 & 0\\
                 0 & \lambda_2 & 0\\
                 0 & 0  & 0\\
               \end{array}
             \right),\left(
               \begin{array}{ccc}
                 \lambda & 1 & 0\\
                 0 & \lambda & 0\\
                 0 & 0  & 0\\
               \end{array}
             \right)$ or $\left(
               \begin{array}{ccc}
                 \lambda & -\mu  & 0\\
                 \mu & \lambda   & 0 \\
                 0 & 0  & 0\\
               \end{array}
             \right)$.

2. If $G$ does not fix any point in $A^2$, then up to conjugation $X$ has the form of

$\left(
  \begin{array}{ccc}
    0 & 0 & a \\
    0 & \lambda  & 0 \\
    0 & 0 & 0 \\
  \end{array}
\right)$ or $\left(
  \begin{array}{ccc}
    0 & 1 & 0 \\
    0 & 0  & 1 \\
    0 & 0 & 0 \\
  \end{array}
\right)$
\end{prop}

\noindent {\bf Proof: }If $G$ fixes a point $P\in A^2$, then $Exp(tX)P=P, \forall t\in \mathbb{R}$. By choosing $P$ as origin, we have $Exp(tX)$ fix the origin. So $X$ must have the form listed in 1.

If $G$ does not fix any point in $A^2$, then %for all $t\in \R$, $r=exp(tX)r+r_0$ has no solution. Hence
$X$ will have an eigenvalue $0$. Computation shows that $X$ must have the form listed in 2.

\begin{coro}The one parameter subgroup corresponding to $X$ is given by

1. $\left(
               \begin{array}{ccc}
                 e^{\lambda_1 t} & 0 & 0\\
                 0 & e^{\lambda_2 t} & 0\\
                 0 & 0  & 1\\
               \end{array}
             \right),\left(
               \begin{array}{ccc}
                 e^{\lambda t} & te^{\lambda t} & 0\\
                 0 & e^{\lambda t} & 0\\
                 0 & 0  & 1\\
               \end{array}
             \right)$ or $\left(
               \begin{array}{ccc}
                 e^{\lambda t}\cos(\mu t)& -e^{\lambda t}\sin(\mu t)  & 0\\
                 e^{\lambda t}\sin(\mu t) & e^{\lambda t}\cos(\mu t)  & 0 \\
                 0 & 0  & 1\\
               \end{array}
             \right)$.

2. $\left(
               \begin{array}{ccc}
                 1 & 0 & at\\
                 0 & e^{\lambda t} & a(e^{\lambda t}-1)\\
                 0 & 0  & 1\\
               \end{array}
             \right)$ or $\left(
               \begin{array}{ccc}
                 1 & t & \frac{1}{2}t^2\\
                 0 & 1 & t\\
                 0 & 0  & 1\\
               \end{array}
             \right)$
\end{coro}
Given a point $P(x_0,y_0)\in A^2$, the curve $Exp(tX)P$ has constant curvature.

\begin{theorem}
Up to affine congruence, the curves generated by the action of one parameter subgroups of $Aff(2)$ on $A^2$ have the forms

1. $y=x^a,a\not=0$. %,1/2,1,2.$

2. $y=ax+bx \ln|x|$;

3. $\displaystyle{\sqrt{(x^2+y^2)}=e^{\frac{a}{b}\arctan(\frac{y}{x})}}$. In polar coordinates it is $\displaystyle{r=e^{\frac{a}{b}\theta}}$. If $a=0$, then this is a circle.
%For differential equation $\displaystyle{\frac{dy}{dx}=\frac{bx+ay}{ax-by}}$, the integral curve at $(1,0)$ is

4. $y=e^x$;\label{curves}
\end{theorem}
%In this section we give a classification of constant curvature curves in $A^2$.

\begin{remark}
The curve $y=x^a$ is affine congruence to the curve $y=x^{\frac{1}{a}}$, and the curve $y=ax+bx \ln|x|$ is affine congruence to the curve $y=x\ln(x)$; The first three cases
of Theorem 6.1 correspond to the three cases in 1 of Corollary 6.1. The fourth case of Theorem 6.1 corresponds to the first case in 2 of Corollary 6.1 with $\lambda\not=0$. The case with $\lambda=0$ and
the second case in 2 of Corollary 6.1 give line and parabola respectively, and they are contained in the first case in Theorem 6.1.
\end{remark}
%
%If a curve $C$ has the constant curvature $k$, then by the moving equation of left moving frames, we have

%$$\displaystyle{\frac{d}{ds} (e_1,e_2)=(e_1,e_2) \left(\begin{array}{cc}k & \frac{1}{3} \\  -1 & \frac{1}{2}k \\   \end{array}   \right)}$$

%This give a differential equation $\frac{dB}{ds}=BK$, choose coordinates such that $B(0)=I$, then we have $B(s)=\exp(tB)$. The Jordan standard form of $B$ is

\begin{theorem}
If $C$ has constant curvature, then $C$ is affine congruence to one of the curve in Theorem \ref{curves}.
\end{theorem}

\noindent{\bf Proof: }We assume $C$ is a smooth curve with constant curvature $k$, in fact $C$ has continuous derivative of order $6$ is enough.
%We have the Cartan structure equation $\D{\frac{d}{ds}(t(s),n(s))=(t(s),n(s))K}$.
Let $F=\left(
                                                                                           \begin{array}{ccc}
                                                                                             t(s) & n(s) & r(s) \\
                                                                                             0 & 0 & 1 \\
                                                                                           \end{array}
                                                                                         \right)$. By Theorem 5.1 the moving equation of $C$ can be written as
$\D{\frac{dF(s)}{ds}=FX}$, where $X=\left(
                       \begin{array}{ccc}
                         0 & -\frac{1}{2} k^2-\frac{1}{2}\frac{dk}{ds}-\frac{\sigma}{3} &  1\\
                         1 & -\frac{3}{2}k& 0\\   0&0&0\\
                       \end{array}
                     \right)$.
%hence $t(s)=(t(0),n(0))Exp(sK)\left(\begin{array}{c} 1 \\   0 \\     \end{array}    \right)$.
Its solution is $F(s)=F(0)Exp(sX)$. By taking the third column of $F$ we get $\left(
                                                                           \begin{array}{c}
                                                                             r(s) \\
                                                                             1 \\
                                                                           \end{array}
                                                                         \right)=F(s)\left(
                                                                           \begin{array}{c}
                                                                             0 \\ 0 \\ 1 \\
                                                                           \end{array}
                                                                         \right)=F(0)Exp(sX)\left(
                                                                           \begin{array}{c}
                                                                             0 \\ 0 \\ 1 \\
                                                                           \end{array}
                                                                         \right)$.
Hence $C$ is affine congruence to a curve listed in Theorem \ref{curves}.

\begin{theorem}

The invariants of the curves with constant curvature are given in the following.

1. For $y=x^a$, $\displaystyle{k(x)=\frac{2\sqrt{3}}{3}sgn(xa (a-1)(a-\frac{1}{2})(a-2)) \frac{(a+1)}{\sqrt{|(2a-1)(a-2)|}}}$. For $a=0$,$\frac{1}{2},1,2$, $S_2=0$ and $k$ is not well defined.
If $\frac{1}{2}<a<2$ then $\sigma>0$, otherwise $\sigma<0$.

%1. For $y=x^a$, $S_2={\frac { \left( {x}^{a} \right) ^{2}{a}^{2} \left( 2\,a-1 \right)
% \left( a-2 \right)  \left( a-1 \right) ^{2}}{{x}^{6}}}$. So

2. For $y=bx+ax\ln|x|$, $\displaystyle{k(x)=-\frac{4\sqrt{3}}{3}sgn(ax)}$ and
$\sigma>0$.

3. For $\displaystyle{\sqrt{(x^2+y^2)}=e^{\frac{a}{b}\arctan(\frac{y}{x})}}$, $\displaystyle{k(x)=-\frac{4\sqrt{3}}{3}\frac{a sgn(b)}{(a^2+9b^2)^{\frac{1}{2}}}}$ and $\sigma>0$.

4. For $y=e^x,\displaystyle{k(x)=\frac{\sqrt{6}}{3}}$ and $\sigma<0$.

\end{theorem}

\noindent{\bf Proof: } Direct computation.
%The characteristic polynomial of $K$ is $\D{\lambda^2-\frac{3}{2} k \lambda +(\frac{1}{2}k^2+\frac{1}{3})=0}$. Its discriminant is $\D{\Delta=\frac{1}{4}(k^2-\frac{16}{3})}$. If $\D{|k|=\pm \frac{4\sqrt{3}}{3}}$,
%then $\Delta=0$. The three cases $\Delta>,=,<0$ correspond to the three cases listed in Theorem ??.

\begin{coro}
If $C$ is a regular curves with zero curvature, then $C$ is part of an eclipse or a hyperbola.
\end{coro}

\section{Modular invariants}
We have defined expression $S_1,S_2$ before. Let $S_3=-45y_{xx}y_{xxx}y_{xxxx}+9y_{xx}^2 y_{xxxxx}+40y_{xxx}^3$, then direct computation shows

$S'_1=\Delta \Gamma^{-3}S_1$;

$S'_2=\Delta^2 \Gamma^{-8}S_2$;

$S'_3=\Delta^3 \Gamma^{-12}S_3$.

We give the following definition.

\begin{define}
Let $f(y,y_x,y_{xx},\cdots)$ be a rational expression of $y(x)$ and its derivatives with $x$, we call $f$ a modular invariant of weight $(p,q)$ if under any affine transformation $f(y',y'_{x'},y'_{x'x'},\cdots)=\Delta^p \Gamma^{q} f(y,y_x,y_{xx},\cdots)$.
\end{define}
%$\displaystyle{y'_{x'x'}=\Delta y_{xx}\Gamma^{-3}}$

%$3y'_{x'x'}y'_{x'x'x'x'}-5y'^2_{x'x'x'}=\Delta^2\Gamma^{-8}(3y_{xx}y_{xxxx}-5y^2_{xxx})$

%$-45y_{xx}y_{xxx}y_{xxxx}+9y_{xx}^2 y_{xxxxx}+40y_{xxx}^3=\Delta^3\Gamma^{-12}(-45y_{xx}y_{xxx}y_{xxxx}+9y_{xx}^2 y_{xxxxx}+40y_{xxx}^3)$

By this definition, $S_1,S_2,S_3$ are modular invariants of weight $(1,-3),(2,-8),(3,-12)$. It is obvious that a modular invariant of weight $(0,0)$ is a differential invariant.
All the modular invariants of weight $(p,q)$ form a real vector space $M_{p,q}$. %which is invariant under the arc derivatives.
$M=\bigoplus\limits_{p,q\in \mathbb{Z}}M_{p,q}$ is a double graded real algebra.

The following proposition determine the structure of $M_{p,q}$.

\begin{prop}
1. $\D{T_1=\frac{S^4_1}{S_3}}$ is a modular invariant of weight $(1,0)$ and $\D{T_2=\frac{S_1 S_2}{S_3}}$ is a modular invariant of weight $(0,1)$.
2. $M_{p,q}=T_1^p T_2^q M_{0,0}$.
\end{prop}

For an element $f\in M_{p,q}$, $f=0$ gives an affine invariant differential operator.
%\section{Affine normal vector and its geometric meaning}

%\section{Orientation}
%If a For an modular invariants $The
%\section{Conclusion}
%In this paper we consider the general affine transformation
%This is our first paper on the general affine geometry. We will continue to study the general affine geometry of curves in higher dimension and study the hypersurfaces in $A^n$. Even the general submanifolds in $A^n$.

%\section{Classification of one parameter subgroups of $Aff(2)$}

\end{document}